\documentclass[12pt]{article} 
\usepackage{amssymb, latexsym,color}
 
\usepackage{times,cite,url}
\newtheorem{defin}{} 
\newtheorem{saetze}[defin]{}
\newtheorem{conjec}[defin]{} 
\newtheorem{lemmas}[defin]{}
\newtheorem{folger}[defin]{} 
\newtheorem{bemerk}[defin]{}

\newcommand{\fillbox}{\mbox{$\bullet$}} 
\newcommand{\ra}{\rightarrow}

\newcommand{\N}{\mathbb N}

\newcommand{\Z}{\mathbb Z}
\newcommand{\Q}{\mathbb Q} 
\newcommand{\G}{\mathcal G}
\newcommand{\B}{\mathcal B} 
\newcommand{\T}{\mathcal T}
\newcommand{\D}{\mathcal D} 
\renewcommand{\O}{\mathcal O}

\newenvironment{items}{\begin{list}{$\alph{item})$} {\labelwidth18pt
      \leftmargin18pt \topsep3pt \itemsep1pt \parsep0pt}} {\end{list}}

\begin{document}

\title{The origins of Coclass Theory} 
\author{Heiko Dietrich and Bettina Eick} 
\date{\today} 
\maketitle

\begin{abstract}
\noindent In 1980, Leedham-Green and Newman introduced the invariant
          {\em coclass} to the theory of groups of prime-power order
          and they proposed five far-reaching conjectures related to
          it. Their work has initiated a deep and fruitful research
          project in group theory that is still ongoing today. We
          outline the main results of this celebrated article, we
          describe the history leading to it, and we survey some
          highlights of work inspired by these results.
\end{abstract}

\section{Introduction}
 
In 1980, Leedham-Green and Newman published ``{\em Space groups and 
groups of prime-power order.\ I\hspace{0.2ex}}'', see \cite{LGN80}, 
in the Archiv der Mathematik. In this article, they considered a new
group-theoretic invariant for groups of prime-power order, called {\em
coclass}, and they formulated five conjectures associated to the
groups of a fixed coclass. These conjectures have inspired a wide
range of interesting research. Remarkably, not only did all
conjectures turn out to be true, their proofs also led to a new
structure theory of $p$-groups and a new approach towards a
classification up to isomorphism. It is fair to say that the Coclass
Article was a turning-point for the theory of finite $p$-groups.
 
We have been invited to review this celebrated article, with the aim  
to briefly explain its main results, to describe the developments
that led to the Coclass Conjectures, and to outline the ongoing
impact on $p$-group theory. Many of the historic details mentioned
below are from an interview with Charles Leedham-Green and Mike
Newman that we have held in November~2022.

In the remainder of this section, we recall the main results of the
Coclass Article.

\subsection{Coclass}\label{sec_cc}
  
The coclass of a group of prime-power order $p^n$ and nilpotency class
$c$ is defined as $n-c$. With the exception of groups of order $p$,
which have coclass $0$, the coclass of a finite $p$-group is an
integer between $1$ and $n-1$. The groups of coclass $1$ are exactly
the groups of {\em maximal class}. These groups have been studied
extensively in the literature (we provide more details below), and
their intricate structure theory was one of the main incentives for
Leedham-Green and Newman to generalise \emph{maximal class} to
\emph{coclass}.  Leedham-Green and Newman also considered this
generalisation as a new approach towards a classification of finite
$p$-groups; they wrote in \cite{LGN80}:

\begin{center}
\begin{minipage}[c]{13cm}
{\em While it seems hopeless to expect a classification up to
  isomorphism [of the $p$-groups of coclass $r$], it may be reasonable
  to hope for a useful classification up to an ``error term''
  depending only on $p$ and $r$.}
\end{minipage}
\end{center} 

\subsection{Conjectures}

The five conjectures proposed by Leedham-Green and Newman in
\cite{LGN80} are nowadays known as the {\em Coclass Conjectures}. They
are decreasing in strength, even though this is not obvious in some
cases.  We recall the conjectures here as follows. Throughout, $p$ is
a prime and $r$ is a positive integer. Conjectures A and B are quoted
from \cite{LGN80}.
\medskip

\noindent {\bf Conjecture A. }  {\it Given $p,r$ there is a positive
integer $f$ such that every $p$-group with coclass $r$ has a normal
subgroup of class $2$ with index at most $p^f$.}
\medskip

\noindent {\bf Conjecture B. }  {\it For every prime $p$ and every
positive integer $r$ there is a positive integer $g$ such that every
$p$-group with coclass at most $r$ has soluble length at most $g$.}
\medskip

Conjectures C, D, and E propose results on the structure of the
\emph{coclass graph} $\G(p,r)$ associated with the finite $p$-groups
of coclass $r$. The vertices of $\G(p,r)$ correspond to the finite
$p$-groups of coclass~$r$, one of each isomorphism type, and the edges
are pairs $(Q, P)$, such that $P$ is an epimorphic image of $Q$ with
$|Q|=p|P|$. Since $Q$ and $P$ have the same coclass, this means that
the classes of $Q$ and $P$ differ by $1$. The definition of edges in
$\G(p,r)$ makes it natural to consider the linear inverse systems
associated with infinite paths in $\G(p,r)$. If $G$ is the inverse
limit of such a system, then $G$ is an infinite pro-$p$-group of
coclass $r$. It is known that this construction induces a one-to-one
correspondence between the maximal infinite paths in $\G(p,r)$ and the
isomorphism types of infinite pro-$p$-groups of coclass $r$. Using
this correspondence, we state Conjectures C, D, and E in the notation
of infinite pro-$p$-groups of coclass $r$.
\medskip

\noindent 
{\bf Conjecture C. }  {\it Every pro-$p$-group of finite coclass is
soluble.}
\medskip

\noindent 
{\bf Conjecture D. }  {\it For fixed $p$ and $r$ there are only
finitely many isomorphism classes of infinite pro-$p$-groups of
coclass $r$.}
\medskip

\noindent 
{\bf Conjecture E. }  {\it For fixed $p$ and $r$ there are only
finitely many isomorphism classes of infinite soluble pro-$p$-groups
of coclass $r$.}
\medskip

It is easy to see that Conjecture A implies Conjecture B, and
Conjecture B implies Conjecture~C. Moreover, if Conjectures C and E
are true, then so is Conjecture D.

\subsection{Space Groups}

The Coclass Article also contains a first investigation of the soluble
infinite pro-$p$-groups of finite coclass, and it exhibits a link to
the theory of space groups. Recall that a space group is an extension
of a translation subgroup $T \cong \Z^d$ by a finite point group $P$
that acts faithfully on~$T$. The theory of space groups has a long
history in group theory with many interesting applications, for
example in crystallography.
A \emph{$p$-adic space group} is a variation of this: it is an
extension of a translation subgroup $T \cong \Z_p^d$, where $\Z_p$
denotes the $p$-adic integers, by a point group $P$ of $p$-power order
that acts faithfully on $T$. Leedham-Green and Newman proved the
following in \cite[Section 4]{LGN80}.
\medskip

\noindent
{\bf Theorem 1.}  {\it Let $G$ be an infinite soluble pro-$p$-group of
  coclass $r$. The hypercentre $H$ of $G$ is finite of order $p^h$
  with $h < r$. The quotient $G/H$ is an infinite soluble
  pro-$p$-group of coclass $r-h$, and it has the structure of a
  $p$-adic space group. }

\section{The history of the Coclass Article}

What insights have led to the proposal of the invariant ``coclass''
and the formulation of the Coclass Conjectures? The following history
combines some collected facts with comments from the personal
interview with Charles Leedham-Green and Mike Newman.

Charles Leedham-Green received his DPhil in 1966 in Oxford and moved
to Queen Mary College (London, England) shortly thereafter, where he
started a collaboration with his also newly arrived colleague Sue
McKay.  Huppert's first book on finite groups \cite{Hup67} had just
appeared, and they read about Blackburn's and Alperin's results on
$p$-groups of maximal class \cite{Alp62,Bla58}. This inspired
Leedham-Green and McKay to work on $p$-groups of maximal class, which
resulted in a sequence of papers \cite{LMc76, LMc78, LMc78a, LGM4}
published between 1976 and 1984.

Mike Newman completed his PhD in 1960 in Manchester.  After his PhD,
Newman took up a position at the Australian National University
(Canberra, Australia).  His interests included $p$-groups and later
also computations.  At the Australian National University he was
involved in the development of the $p$-group generation algorithm,
see \cite{New77} for an early outline. Computational experiments with
this algorithm were very important for many advances in Coclass Theory
(and beyond).

\subsection{Writing the Coclass Article: 1974--1980}

Leedham-Green and Newman met at the ICM 1974 in Vancouver where both
of them gave a `short communication': Newman talked about {\em Groups
  of exponent four} and Leedham-Green reported {\em On $p$-groups of
  maximal class}. As both shared an interest in $p$-group theory, they
naturally started a collaboration. In his 1990 paper
\cite[p.\ 52]{New90}, Newman recalled that the lecture of
Leedham-Green was the first time he has seen a coclass graph for
groups of maximal class. At that time, it was known that the graph
$\G(p,1)$ has a unique maximal infinite path, and the inverse limit of
the groups on this path yields the infinite pro-$p$-group $G_p$ of
maximal class. Specifically, $G_p$ is isomorphic to a semidirect
product $P\ltimes T$ where $P$ is a cyclic group of order $p$ acting
via multiplication by a $p$-th root of unity on the
$(p-1)$-dimensional $p$-adic lattice $T$.

In 1976, Leedham-Green visited Newman in Canberra. On the way out to
Canberra, it occurred to Leedham-Green that the infinite pro-$p$-group
associated with $\G(p,1)$ can be generalised by replacing $P$ by a
cyclic group of order $p^n$, acting via multiplication by a $p^n$-th
root of unity on a $(p-1)p^{n-1}$-dimensional $p$-adic lattice. Upon
arrival in Canberra, Leedham-Green discussed this generalisation with
Newman who recognised this type of construction as a $p$-adic
variation of a space group. Newman had just returned home from a visit
to the RWTH Aachen, where Neub\"user and his collaborators Brown,
B\"ulow, Zassenhaus, and the crystallographer Wondratschek were
writing a book on the classification of space groups in dimension~$4$.
Newman transferred many results for integral space groups to $p$-adic
space groups. Similar to how $G_p$ yields an infinite path in
$\G(p,1)$, the newly constructed infinite groups gave rise to infinite
families of finite $p$-groups, all having the same coclass. This
combination of insights gave birth to the new invariant!
  
Conjectures A and B generalise results for $p$-groups of maximal class
and have been influenced by the work of Wiman \cite{Wim46, Wim52},
Blackburn \cite{Bla58}, Alperin \cite{Alp62}, Miech
\cite{Mie70,Mie74a,Mie78,Mie82}, Shepherd \cite{She71}, and
Leedham-Green and McKay \cite{LMc76, LMc78, LMc78a,
  LGM4}. Particularly relevant to Conjecture B is a theorem by Alperin
that states: {\em For every $p$, there is a constant $m$ such that
  every $p$-group of maximal class has derived length at most $m$}.

Conjectures C, D, and E are of a different nature and rely on the new
connection between $p$-groups and space groups. Leedham-Green wrote in
\cite{LGM94a}:

\begin{center}
\begin{minipage}[c]{13.2cm}
{\em These conjectures did not come out of the blue.  They were first
  formulated (inaccurately) from the insights obtained from my work
  with Susan McKay on $p$-groups of maximal class (that is, of coclass
  1), which gave us a proof of the conjectures in this case. That the
  conjectures appeared in the correct form (and I am still surprised
  at how accurate they were) was due to M.~F.~Newman's insight into
  the structure of space-groups.}
\end{minipage}
\end{center}

\subsection{Immediate Reactions: 1975--1984}

The results and conjectures of the Coclass Article had been talked
about prior to the submission in 1979. Many researchers were
eager to join the project, and some had already proved some
interesting results. It was a time similar to the {\em gold-rush}, and
it must have been extraordinarily exciting.

The suffix ``I'' of the Coclass Article suggests that the paper was
meant to be the first of a series.  Indeed, immediately after its
publication, a sequence of papers appeared in print.  The first two of
them (called ``II'' and ``IV'') where written at the RWTH Aachen
\cite{spaceII, spaceIV}. The first by Finken, Neub\"user, and Plesken
provides a solution to the isomorphism problem of the $p$-adic space
groups with fixed coclass. The second by Felsch, Neub\"user, and
Plesken determines an infinite family of counterexamples to the
class-breadth conjecture via quotients of $p$-adic space groups. Part
``III'' of the series is mentioned in \cite{spaceII}, but was never
published. The other articles of the series (called ``V'', ``VI'',
etc.) were aimed at proving the Coclass Conjectures; we provide more
details below.

Very little was known beyond the detailed structure theory of the
$p$-groups of maximal class and thus more explicit investigations were
of interest.  In 1975, James \cite{James, James2} started the
investigation of the groups in $\G(2,2)$.  He employed algorithmic
methods and even early computers to solve the isomorphism problem for
these $2$-groups. This was one of the starting points for Ascione
\cite{ascione, ascione2, ascionephd}, a PhD student of Newman, to
consider the $3$-groups of coclass~$2$. Her thesis (137 pages long)
provided interesting insights into $\G(3,2)$.

\section{The proof of the Coclass Conjectures}

The proofs of the Coclass Conjectures emerged in the years
1980--1994. In many cases, the proofs or at least significant ideas for
proofs were around years before they were finally published. Hence, the
publication dates often do not accurately reflect the historic
events. A full account of the proofs is available in the book by
Leedham-Green and McKay \cite{LGM02}.

\subsection{The proof of Conjecture E}   

At the time when the Coclass Article was written, a proof of
Conjecture E was already in the air. It follows from Theorem 1 that
Conjecture E is equivalent to the claim that there are only finitely
many isomorphism classes of soluble $p$-adic space groups of fixed
coclass. Earlier, around 1911, Bieberbach had proved that there are
only finitely many isomorphism classes of integral space groups in any
fixed dimension. Thus, it remained to study the relationship between
soluble $p$-adic space groups and integral space groups, and to
describe how the coclass of the group is related to the dimension of
the space group.

Shortly after the Coclass Article appeared, Leedham-Green, McKay, and
Plesken met at the Mathematisches Forschungsinstitut in Oberwolfach to
discuss ideas for a proof. One step in the proof was the determination
of all the possible point groups; that is, finite $p$-groups $P$ that
can be embedded into ${\rm GL}(d, \Z_p)$ for some suitable dimension
$d$. A crucial ingredient for this step was Vol'va\v{c}ev's
description \cite{volv1} of the Sylow $p$-subgroups of ${\rm GL}(d,
\Q_p)$. Leedham-Green and Plesken \cite{LGP} corrected an error for
$p=2$, and eventually completed the determination of the possible
point groups.

The investigation of the extensions of translation subgroups by point
groups and the relation between dimension and coclass turned out to be
a more subtle business. Charles Leedham-Green remembers:

\begin{center}
\begin{minipage}[c]{13cm}
{\em Some while later I found myself in Aachen, working with Wilhelm
  Plesken on Conjecture E. We formulated a conjecture that would
  settle the matter, and tested it in dimension 4, using a computer
  program commissioned by the Berlin bus company.  The buses ran on
  hydrogen that was stored in rare earth crystals, and the hydrogen
  reduced their symmetry.  So a program was commissioned to compute
  subgroups of space groups.  Using it, the machine found a
  counter-example to our conjecture. }
\end{minipage}
\end{center}

Despite such problems, Leedham-Green, McKay, and Plesken finally
obtained a proof of Conjecture E.  The proof appeared in two papers
\cite{spaceV, spaceVI}, distinguishing between even and odd
primes. For odd primes, there is only one Sylow $p$-group in 
${\rm GL}(d, \Q_p)$, and the extensions constructed from its subgroups
behave in a rather generic way. For the even prime $p=2$, this is no
longer true and the proof required more work. Explicit bounds on the
coclass for such $p$-adic space groups have later been proved by McKay
\cite{McKay94}.

\subsection{The proof of Conjecture C}

The next step was the proof of Conjecture C; together with the proof
of Conjecture E, this also settled Conjecture D.

A first milestone was Donkin's proof \cite{spaceVIII} of Conjecture C
for primes $p \geq 5$.  This proof relies on the classification of
simple $p$-adic Lie algebras and the theory of $p$-adic Chevalley
groups. Donkin \cite{spaceVIII} stated that a proof covering {\em all}
primes has been suggested by Tits; this proof uses the theory of
buildings, but it remained unpublished.

The second milestone was the comparatively short and significantly
more elementary proof by Shalev and Zelmanov \cite{SZ92} holding for
all primes.  It relies on a reduction to analytic groups obtained by
Leedham-Green \cite{LGM94a}.

\subsection{The proof of Conjecture A}

Two independent proofs for Conjecture A have been found by
Leedham-Green \cite{LGM94b} and Shalev \cite{Sha94}. Both proofs
emerged around 1990/91 and appeared in print in 1994. Mann has written
an excellent `featured review' of both proofs on the AMS MathSciNet.
To underline the significance of Conjecture A, Mann also states:

\begin{center}
\begin{minipage}[c]{13cm}
{\em According to a reported saying of N.\ Blackburn, this conjecture,
  once proved, will be 'the first general theorem of the theory of
  $p$-groups'.}
\end{minipage}
\end{center}

Prior to \cite{LGM94b,Sha94}, special cases of Conjecture A had been
established. For example, McKay \cite{McKay87} proved Conjecture A for
so-called uncovered CF-groups. Based on this proof, and using the
newly developed theory of powerful $p$-groups, Mann \cite{spaceVII}
proved Conjecture A for all uncovered $p$-groups if $p>2$.

Shalev's proof \cite{Sha94} of Conjecture A also uses the theory of
powerful $p$-groups and relies on the results of Mann
\cite{spaceVII}. More precisely, it uses that every finite group of
coclass $r$ has a powerful subgroup of index bounded by $p$ and
$r$. Assuming Conjecture A is false, these powerful subgroups in a
family of assumed counterexamples are used to produce an appropriate
Lie ring that, eventually, contradicts a result of Jacobson on Engel
Lie algebras.

Leedham-Green's proof \cite{LGM94b} appeared later than Shalev's, but
(according to Mann) it was the first proof and it was available to
Shalev. It relies on Conjecture C and its proof. Inverse limit
constructions are a central step in the proof: Assuming Conjecture~A
is false, an inverse limit construction of an assumed family of
counterexamples produces an infinite pro-$p$-group $G$ of coclass
$r$. This is $p$-adic analytic by \cite{LGM94a}, and thus Conjecture C
and its proof have implications on the structure of $G$.  The final
contradiction is reached by employing results for so-called settled
groups and the Lie rings associated with them.

\section{Beyond the Coclass Conjectures}

After the Coclass Conjectures have been proved, the main focus shifted
on the detailed structure of the coclass graphs $\G(p,r)$. The proof
of Conjecture D implies that there are only finitely many maximal
infinite paths in $\G(p,r)$ for any given $p$ and $r$. For a group $G$
in $\G(p,r)$ denote by $\D(G)$ the full subtree of $\G(p,r)$
consisting of $G$ and its descendants. If $G$ is a group on an
infinite path, then $\D(G)$ is an infinite tree. It is called a {\em
  coclass tree} if it contains only one infinite path starting at its
root $G$. It is a {\em maximal coclass tree} if there is no proper
ancestor $H$ of $G$ in $\G(p,r)$ such that $\D(H)$ is also a coclass
tree. The proof of Conjecture D now implies the following; see also
the discussion in \cite[Section~11.2]{LGM02}:
\medskip

\noindent {\bf Theorem D'.} {\it The graph $\G(p,r)$ consists of
  finitely many maximal coclass trees and finitely many other groups.}
\medskip

Theorem D' implies that the general structure of a coclass graph is
determined by the structure of its finitely many maximal coclass
trees. By definition, every coclass tree $\T$ has a unique infinite
path, $S_1 \ra S_2 \ra \ldots$ say, starting at its root. The inverse
limit $S$ of the groups on this path is the infinite pro-$p$-group of
coclass $r$ associated with this maximal coclass tree. For $n \in \N$
we define the {\em branch} $\B_n$ of $\T$ as the difference graph 
$\D(S_n) \setminus
\D(S_{n+1})$; that is, $\B_n$ consists of all descendants of $S_n$ in
$\G(p,r)$ that are not descendants of $S_{n+1}$.

For odd $p$ and any given $r$, Eick \cite{Eic04b} introduced an
algorithm to determine, up to isomorphism, the infinite pro-$p$-groups
of coclass $r$; this can be considered as a constructive version of
the proof of Theorem~D.  Further, Eick \cite{Eic05} determined
estimates for the number of isomorphism classes of infinite
pro-$p$-groups of coclass $r$, and showed that this number grows
exponentially with $p$ and~$r$.

Regarding the structure of a coclass tree, the following observation
holds.

\medskip

\noindent
{\bf Remark 1.}  {\it Each maximal coclass tree $\T$ consists of its
  maximal infinite path $S_1 \ra S_2 \ra \ldots$ and its branches
  $\B_1, \B_2, \ldots$.  Each branch $\B_i$ is a finite tree with root
  $S_i$.}
\medskip

By Theorem D' and Remark 1, the general structure of a coclass graph
$\G(p,r)$ is determined by the structure of the branches of the
maximal coclass trees in $\G(p,r)$. Understanding these branches is
the main focus in coclass theory today.

\subsection{Skeletons} 
\label{sec_const}

In their papers on $p$-groups of maximal class, Leedham-Green and
McKay introduced a novel construction for certain groups in the
branches of the coclass tree in $\G(p,1)$: the so-called {\em
  constructible groups}.  This idea was later generalised to all
graphs $\G(p,r)$ as described in the book \cite{LGM02}. The
constructible groups contained in a branch $\mathcal{B}_i$ form a
connected subtree. This subtree, truncated at a depth close to the full
depth,  is called the {\em skeleton} of $\mathcal{B}_i$. The following
theorem of Leedham-Green shows that the skeletons exhibit the broad
structure of the branches, see \cite[Theorems 11.3.7 and 11.3.9]{LGM02} 
for a proof. This can be considered as a proof for a stronger version 
of  Conjecture A.
\medskip

\noindent {\bf Theorem A'.}  {\it Given $p,r$, there are integers $n =
  n(p,r)$ and $m = m(p,r)$ such that every $p$-group $P$ of coclass
  $r$ of order at least $p^n$ has a normal subgroup $N$ of order at
  most $p^m$ so that $P/N$ is a constructible group.}
\medskip 

The explicit definition of constructible groups is technical and goes
beyond the scope of this survey. The underlying idea is to consider
the groups on the infinite path of a coclass tree. These groups are
quotients of the associated inverse limit $S$, and the infinite
pro-$p$-group $S$ is an extension of a normal subgroup $T \cong
\Z_p^d$ by a finite $p$-group $P$. The group $T$ has a unique maximal
$S$-invariant chain $T = T_0 > T_1 > \ldots$, and the groups on the
original infinite path are isomorphic to the quotients $S/T_i$ for
large enough $i$. The constructible groups are now obtained by
`twisting' the multiplication in $T/T_i$ such that the result is a
group of class $2$. The `twisting' into a group of class $2$ is
facilitated via the surjective homomorphisms in
\[ {\rm Hom}_{S/T}(T \wedge T, T_j/T_i).\]
The choice of $i$ prescribes the order of the resulting constructible
group, whereas the choice of $j$ prescribes the branch in which the
group lies. A detailed study of the homomorphism space and the
isomorphism problem for constructible groups translates to interesting
and non-trivial problems in algebraic number theory; this has been
investigated in many publications, most recently in
\cite{de17,DS19,CDEM21}.

Theorem A' implies that a classification of skeleton groups yields a
classification of all $p$-groups of coclass $r$, up to a small
``error-term''; it provides a solution for one of the main aims in
Coclass Theory, see the quote in Section \ref{sec_cc}.

\subsection{The graph $\G(5,1)$ and computations}
\label{sec_comp}

Computational group theory had a major influence in Coclass
Theory. For example, finite parts of coclass graphs can be computed
with the $p$-group generation algorithm, which was first developed by
Newman \cite{New77} and later significantly extended by O'Brien
\cite{OB90}. Computer experiments with the $p$-group generation
algorithm were crucial for formulating many of the most recent
conjectures on the structure of $\G(p,r)$.

The graph $\G(2,1)$ is illustrated in \cite{LGN80} based on the
well-known result that the number of $2$-groups of maximal class of
order $2^n$ is $3$ provided that $n\geq 4$. The graph $\G(3,1)$ can
be read off from the determination by Blackburn \cite{Bla58} of the $3$-groups
of maximal class. For larger $p$, Leedham-Green and McKay \cite{LMc76,
  LMc78a, LMc78, LGM4} described many features of $\G(p,1)$. However,
extensive details on the precise structure of $\G(p,1)$ were not known
around the year 1990. At that time, computational methods were
employed to investigate the smallest open case, the graph
$\G(5,1)$. Newman's experiments \cite{New90} revealed that $\G(5,1)$
is significantly more complex than $\G(2,1)$ and $\G(3,1)$, yet it has
a lot of periodic structure. To describe the latter, we need the
notion of \emph{pruned branches}: for a branch $\B_i$ and an integer
$k$, the pruned branch $\B_{i}(k)$ is the full subtree of $\B_i$ of
all vertices of distance at most $k$ from the root. The following is a
simplified variation of \cite[Conjecture~IV]{DEF08}.

\pagebreak

\noindent
{\bf Conjecture ($\G(5,1)$).}  {\it If $i\geq 9$, then the branches
  $\B_i$ of the unique maximal coclass tree $\T$ in $\G(5,1)$ have
  depth $i+1$, and the following hold:
\begin{items}
\item[\rm (a)] $\B_{i+4}(i-4) \cong \B_i(i-4)$, and
\item[\rm (b)] $\B_{i+4} \setminus \B_{i+4}(i-4) \cong \B_i \setminus
  \B_i(i-8)$.
\end{items}}

If this conjecture is true, then $\G(5,1)$ can be constructed from a
finite subtree using two types of periodicity patterns, both with
period $4$. The first pattern is part (a): it moves parallel to the
infinite path in $\G(5,1)$ and describes how a large pruned part of
$\B_i$ can be embedded into $\B_{i+4}$. The second pattern is (b): it
moves vertically to the infinite path and describes how the ``feet''
of the branches growth.

Part (a) and a weaker form of part (b) have been proved by Dietrich
\cite{hdmc1,hdmc2}, while the complete proof is still open today.

\subsection{The graphs $\G(2,r)$ and Conjecture P} 

As often in mathematics, the case $p=2$ behaves differently to odd
primes. The most remarkable difference is that all the skeletons of a
maximal coclass tree in $\G(2,r)$ are trivial: they only contain a
single vertex, the root of the associated branch. Theorem A' now
implies that the branches in $\G(2,r)$ all have {\em bounded depth},
that is, for every coclass $r$, there is a constant $c$ such that
every branch in $\G(2,r)$ has depth at most $c$. This is not true in
general for coclass trees with $p>2$; for example, see the conjecture
on $\G(5,1)$ above. This shows that the graphs $\G(2,r)$ have a much
simpler structure than most of the graphs $\G(p,r)$ for odd primes.

Newman and O'Brien \cite{NOB99} used computational tools to
investigate the coclass trees in $\G(2,r)$ for $r \leq 3$. They
classified the groups in $\G(2,3)$ and proposed the following
conjecture.
\medskip

\noindent 
{\bf Conjecture P. } {\it Let $\T$ be a coclass tree in $\G(2,r)$ with
  branches $\B_1,\B_2,\ldots$.  There exist integers $d$ and $n$ such
  that $\B_i \cong \B_{i+d}$ for each $i \geq n$.}
\medskip

This conjecture was proved almost immediately by du Sautoy
\cite{DuS01} in the following generalised form. The \emph{dimension}
of a coclass tree is the dimension of the $p$-adic space group associated
to the infinite path in the tree.
\medskip
 
\noindent 
{\bf Theorem P'. }  {\it Let $p$ be an arbitrary prime, let $r,k$ be
  positive integers, and let $\T$ be a coclass tree in $\G(p,r)$ with
  branches $\B_1,\B_2,\ldots$ and dimension $d$. There exists an
  integer $n = n(p,r,k)$ such that $\B_i(k) \cong \B_{i+d}(k)$ for all
  $i \geq n$.}
\medskip

Du Sautoy's proof for Theorem P' is based on model theory and the
theory of zeta-functions; it is non-constructive and it does not yield
any bounds for $n$. A second proof for Theorem P' was given by Eick
and Leedham-Green \cite{ELG08}; their proof is constructive and has
the additional feature that it yields an explicit analysis of the
structure of the groups in $\B_i(k)$, as well as an explicit bound for
$n = n(p,r,k)$.

\subsection{The graph $\G(3,2)$ and Conjecture W}

The $3$-groups of coclass $2$ were already investigated by Ascione,
Havas, and Leedham-Green \cite{ascione, ascione2} around the time when 
the Coclass Article was written. Additionally, Leedham-Green wrote 
detailed notes on $\G(3,2)$ which were intended to be Part ``III'' in 
the series of papers started by the Coclass Article. These notes 
were revisited in 2001 by Eick, Leedham-Green, Newman, and O'Brien.  
Their paper \cite{ELNO07} contains the following result.
\medskip
 
\noindent {\bf Theorem ($\G(3,2)$).} {\it
\begin{items}
\item[\rm (a)] The graph $\G(3,2)$ consists of $16$ maximal coclass
  trees. Among the associated infinite pro-$3$-groups, there are six
  that are central extensions of the infinite pro-$3$-group of
  coclass~$1$ and there are ten that are $3$-adic space groups.
\item[\rm (b)] Among the $16$ maximal coclass trees, $12$ have bounded
  depth; their structure is fully described by Theorem P'. The
  remaining four maximal coclass trees correspond to $3$-adic space
  groups whose point group is cyclic of order $9$ or extraspecial of
  order $27$.
\end{items}}
\medskip

The main part of \cite{ELNO07} was concerned with a detailed
description of the four maximal coclass trees that were not of bounded
depth. These descriptions all rely on the full determination of the
skeletons of the branches. Based on this investigation of $\G(3,2)$,
the available experimental evidence for $\G(5,1)$, and the periodicity
results in \cite{hdmc2}, the following conjecture was proposed
\cite[Conjecture W]{ELNO07}:
\medskip

\noindent {\bf Conjecture W. } {\it Let $\T$ be a coclass tree in
  $\G(p,r)$ with branches $\B_1, \B_2, \ldots$ and dimension
  $d$. There are integers $n,k$ such that for all $i \geq n$ the
  following hold:
\begin{items}
\item[\rm (a)] $\B_i(k)\cong \B_{i+d}(k)$, and
\item[\rm (b)] for each $P$ of depth $k$ in $\B_{i+d}$ there exists
  $Q$ of depth $k-d$ in $\B_i$ such that $\D(P) \cong \D(Q)$.
\end{items}}
\medskip

Conjecture W is still open, and proving it is one of the main aims of
coclass theory. The main problem is to find a generic description for
$Q$ based on $P$.  The group $Q$ is often called a {\em periodic
  parent} for $P$.

\subsection{The graphs $\G(p,1)$ and Galois trees}

Let $\T_p$ denote the unique maximal coclass tree in $\G(p,1)$. Then
$\T_p$ has branches of bounded depth for $p \in \{2,3\}$, it has
branches of growing depth but bounded width for $p=5$, and it is
growing wildly for $p \geq 7$. The growth for $p \geq 7$ has recently been
investigated by Dietrich and Eick \cite{de17} and by Cant, Dietrich,
Eick, and Moede \cite{CDEM21}.

For a group $G$ in the skeleton of a branch in $\T_p$ we define its
{\em Galois order} $o(G)$ as the $p'$-part of $|{\rm Aut}(G)|$. If $G$
is on the infinite path of $\T_p$, then $o(G) = (p-1)^2$. If $G$ is
not on the infinite path, then $o(G)$ divides $p-1$. A subtree $\O$ of
a skeleton $\mathcal{S}_i$ of a branch $\B_i$ in $\T_p$ is a {\em
  Galois tree} of order $o$ if all groups in $\O$ have Galois order
$o$ and $\O$ is a subtree that is maximal with this property.
Dietrich and Eick \cite[Theorem 1.1]{de17} determined the Galois trees
with Galois order $p-1$ in $\G(p,1)$ as follows.
\medskip

\noindent
{\bf Theorem (Galois trees).}  {\it Let $p \geq 5$ and $\ell =
  (p-3)/2$. Each skeleton in $\G(p,1)$ has $\ell$ Galois trees with
  maximal Galois order $p-1$. Their roots $G_1, \ldots, G_\ell$ have
  depth $1$ in the skeleton, and their leaves have full depth in the
  skeleton. A group of non-maximal depth $m$ in the Galois tree with
  root $G_i$ has either $1$ or $p$ descendants; the latter case occurs
  if and only if $m$ is even with $(m \bmod (p-1)) \notin \{0,p-1-2i\}$.}
\medskip

Cant, Dietrich, Eick and Moede \cite{CDEM21} generalised this result
to Galois trees with arbitrary Galois order.

The main motivation behind investigating Galois trees is indicated in
Dietrich's work \cite{hdmc2} and in \cite{de17}: it is believed that
Galois trees might facilitate a method to determine a periodic parent
$Q$ for a group $P$ (as proposed in Conjecture W). More precisely, the 
hope is that
periodic parents can be obtained by moving $d$ steps up in a Galois
tree of $\B_{i+d}$ and then taking the corresponding group in
$\B_i$. It is known that moving $d$ steps up in a skeleton will, in
general, not provide a periodic parent, see Saha \cite[Section
  6.2.7]{saha}.

\subsection{Coclass families}

The periodicity result by Eick and Leedham-Green \cite{ELG08}
partitions the infinitely many groups in a sequence of pruned branches
$\B_{1}(k), \B_{2}(k),\ldots$ into finitely many \emph{coclass families}.
For this, let $n$ and $d$ be as in Theorem P' and start with a group $G$ 
in $\B_{i}(k)$ with $i\in\{n,n+1,\ldots,n+d-1\}$. Then the associated 
coclass family is defined as the sequence of groups that can be reached 
by taking all the images of $G$ under the graph isomorphisms 
$\B_i(k) \ra \B_{i+d}(k) \ra \B_{i+2d}(k)\ra \ldots$ as defined in the
proof of Theorem P' by Eick and Leedham-Green \cite{ELG08} which also
shows that all the groups in one such family can be described by a single 
\emph{parametrised
  presentation}. Well-known examples of this type of presentation are
the coclass families in $\G(2,1)$: these are the families of dihedral
groups $D_{2^n}$, quaternion groups $Q_{2^n}$ and semidihedral groups 
$SD_{2^n}$ and they can be described by the following parametrised
presentations:
\[
\begin{array}{lcll}
D_{2^{n}} &= &\langle x,y \, | \, x^{2^{n-1}}=1,\;y^2=1,\; x^y=x^{-1}
\rangle&(n\geq 3),\\ Q_{2^{n}} &=& \langle x,y \, | \,
x^{2^{n-1}}=1,\;y^2=x^{2^{n-2}},\;x^y=x^{-1} \rangle&(n\geq
4),\\ SD_{2^{n}} &=& \langle x,y \, | \, x^{2^{n-1}}=1,\;y^2=1,\;
x^y=x^{2^{n-2}-1} \rangle&(n\geq 4).
\end{array}
\]

In turn, the parametrised presentations allow one to prove results for
an infinite family of groups.  Results of this flavour have been
obtained for automorphism groups, Schur multipliers, cohomology rings,
character degrees, and various other invariants, see, for example, 
\cite{carl,Couson14,Eic06b, DEF08, Eick08, schur2, EickSchur, New09, 
cohom, Sym21, EG17, assalg, mor, pres2}.

\subsection{Other algebraic objects}

As final comment, we mention that the invariant ``coclass'' has also
been defined for other algebraic structures, such as nilpotent Lie
algebras, nilpotent associative algebras, and semigroups, see for
example \cite{riley,assalg2, assalg,semi}.  These investigations have
also revealed interesting results, periodicity patterns, and various
conjectures.

\section{Open problems and further research}

We list a few open problems which
we believe are worth to follow up.
Many details and background on these will be available in the 
forthcoming book \cite{book2}.
\medskip

\noindent
{\bf Problem 1.} {\it Problem 17.63 in the Kourovka notebook.}
\\ Prove that if $p$ is an odd prime and $s$ a positive integer, then
there are only finitely many isomorphism classes of $p$-adic space
groups of finite coclass with point group of coclass $s$.
\medskip

\noindent
{\bf Problem 2.} {\it Problem 14.64 in the Kourovka notebook and
  Problem 3 in Shalev's survey \cite{Sha95}.} \\ Complete the
description of $\G(5,1)$ and the classification of the $5$-groups of
coclass~$1$ via coclass families.
\medskip

\noindent
{\bf Problem 3.} {\it More general than Problem 2 is: } \\ Complete
the description of $\G(p,1)$ for $p>5$ and the classification of the
$p$-groups of maximal class via coclass families.
\medskip

\noindent
{\bf Problem 4.} {\it Even more general than Problem 3:} \\ Prove
Conjecture W.


\section*{Acknowledgments}

We thank Charles Leedham-Green and Mike Newman for introducing us to
the exciting area of coclass theory. We also thank them for providing
interesting background information, details, and anecdotes about
historic dates and events that led to their 1980 paper.

\end{document}